\documentclass{amsart}
\usepackage{amssymb}
\usepackage{amsthm}
\newtheorem{thm}{Theorem}[section]
\newtheorem{lem}[thm]{Lemma}
\newtheorem{prop}[thm]{Proposition}
\newtheorem{cor}[thm]{Corollary}

\newtheorem{remark}[thm]{Remark}

\newcommand{\Z}{\mathbb{Z}}
\begin{document}

\title[Normal Restriction in Finite Groups]{Normal Restriction in Finite Groups}

\author{Hung P. Tong Viet}

\address{The School of Mathematics, The University of Birmingham, Edgbaston, Birmingham, ~B15~2TT,
 United Kingdom}

\email{tongviet@for.mat.bham.ac.uk}

\thanks{}


\date{\today}

\dedicatory{}

\commby{}
\begin{abstract}
A subgroup $H$ of a finite group $G$ is called an $NR-$subgroup (Normal Restriction) if, whenever $K\unlhd H,$ then $K^G\cap H=K,$ where $K^G$ is the normal closure of $K$ in $G.$ In this paper, we will prove some sufficient conditions for the solvability of finite groups which possess many $NR$-subgroups. We also prove a criterion for the existence of a normal $p$-complement in finite groups.
\end{abstract}
\maketitle

\section{Introduction}
\label{intro}
Let $G$ be a finite group. A subgroup $H$ of $G$ is called a \emph{$CR$-subgroup} (\emph{Character Restriction}) of $G$ if every complex irreducible character of $H$ is a restriction of some irreducible character of $G.$ It is well known that if $H$ is a $CR$-subgroup of $G$ and $K\unlhd H,$ then $K^G\cap H=K.$ This leads to the following definitions: a triple $(G,H,K)$ is said to be {\em special} in $G$ if $K\unlhd H\leq G$ and $H\cap K^G=K,$ where $K^G$ is the normal closure of $K$ in $G.$ A subgroup $H$ is called an
$NR-$subgroup (Normal Restriction) if, whenever $K\unlhd H,$ then $(G,H,K)$ is special in $G.$ From definition, we see that if $H\leq G$ and $H$ is simple then $H$ is an $NR$-subgroup of $G.$

Li Shirong \cite{shi} called a subgroup $K,$ an $NE$-subgroup if $(G,N_G(K),K)$ is special in $G.$ He showed that
if every minimal subgroup of $G$ is an $NE$-subgroup then $G$ is solvable. Such a group is called a $PE$-group and the structure of
a minimal non-PE-group was investigated. Yangming Li\cite{yan} showed that if every minimal subgroup of
prime order or of order $4$ is an $NE$-subgroup of $G,$ then $G$ is supersolvable. The author also classified non-abelian simple
groups  whose second maximal subgroups are $PE$-groups.

In Tong-Viet\cite{tong}, it is shown that if every maximal subgroup of $G$ is an $NR$-subgroup then $G$ is solvable. This gives a positive answer to a question posed in Berkovich\cite{ber}. We see that in the symmetric group $S_4$ only self-normalizing non-nilpotent maximal subgroups are $NR$-subgroups but this group is still solvable. Our first result is a generalization of Theorem $1.1$ in Tong-Viet\cite{tong}.

\begin{thm}\label{th1} If every non-nilpotent maximal subgroup of $G$ is either normal or $NR$ in $G,$ then $G$ is solvable.
\end{thm}

A maximal subgroup of $G$ is said to be {\em 1-maximal}. For $n\geq 2,$ a subgroup $H$ is called {\em n-maximal} if it is maximal in some $(n-1)$-maximal subgroup of $G.$ In view of Theorem \ref{th1}, it is natural to ask whether or not a group $G$ is solvable if every non-nilpotent $2$-maximal subgroup of $G$ is either subnormal or $NR$ in $G.$ The answer to this question is `No'. The minimal counter-example is the alternating group $A_5.$ It is easy to see that every $2$-maximal subgroup of $A_5$ is nilpotent. In the next theorem, we will show that in fact $A_5$ is the unique non-abelian composition factor of groups satisfying the assumption of the above question.

Denote by $S(G)$ the solvable radical of a group $G,$ that is, a maximal normal solvable subgroup of $G.$
\begin{thm}\label{th2} If  every  non-nilpotent $2$-maximal subgroup of $G$ is either  subnormal or $NR$ in $G,$ then $G/S(G)$ is trivial or isomorphic to  $A_5.$
\end{thm}
As a consequence, we obtain another sufficient condition for the solvability of finite groups as follows:
\begin{cor}\label{cor} If  every  $2$-maximal subgroup of $G$ is either subnormal or $NR$ in $G$ then $G$ is solvable.
\end{cor}
We cannot extend Corollary \ref{cor} further as all $3$-maximal subgroups of $A_5$ are $NR$ in $A_5.$
Recall that a subgroup $H$ is said to have a \emph{normal complement} in $G$ if there exists a normal subgroup $L$ of $G$ such that $G=HL$ and $H\cap L=1.$
In Isaacs\cite{isa}, it is shown that if $N=N_G(P)$ is a $CR$-subgroup, where $P\in Syl_p(G),$ then $N$ has a normal complement in $G.$ This result has been generalized in Berkovich\cite{ber}. In that paper, the author replaces the property $CR$ by $NR$ and still obtains the same conclusion. However we can see that the argument may apply to a broader class of subgroups, say $p$-subgroups instead of $p$-Sylow subgroups.

\begin{prop}\label{nc1} Let $P$ be a $p$-subgroup of $G$ and $N=N_G(P).$ If the triples $(G,N,P)$ and $(G,N,\Phi(P))$ are special in $G$ then $N$ has a normal complement $T$ in $G,$ with $(p,|T|)=1.$
\end{prop}
Applying Lemma \ref{lem1}$(e)$, we get the following corollary:

\begin{cor}\label{nc2} Let $P$ be a $p$-subgroup of $G$ and $N=N_G(P).$ Assume that for any $T\in \{P,\Phi(P)\},$ there exists a subgroup $L$ of $G$ such that $G=NL$ and $N\cap L=T$ then $N$ has a normal complement in $G.$
\end{cor}

Recall that a group $G$ is said to be \emph{$p$-nilpotent} if there exists a normal $p'$-subgroup $N$ of $G$ such that $G=SN$ where $S\in Syl_p(G)$ for some prime $p$ and we call $N$ a \emph{ normal $p$-complement} in $G.$ As an application of Proposition \ref{nc1}, we prove a new criterion  for the existence of a normal  $p$-complement in finite groups.

\begin{thm}\label{th4} Let $P$ be a $p$-subgroup of $G$ and $N=N_G(P).$ Assume that the triples $(G,N,P)$ and $(G,N,\Phi(P))$ are special in $G.$ If $N$ is $p$-nilpotent then $G$ is $p$-nilpotent.
\end{thm}

Observe that the existence of a supersolvable maximal $NR$-subgroup does not imply the supersolvability of a group. The minimal counter-example is $A_4$ with a supersolvable maximal $NR$-subgroup $A_3.$
\begin{thm}\label{th5} Let $H$ be an $NR$-subgroup of prime index in $G.$ If $H$ is supersolvable then $G$ is supersolvable.
\end{thm}

All groups are finite. We adapt the notations in Conway et al.\cite{atlas} for finite simple groups. All other notations for finite groups are standard.
\section{Preliminaries}
\label{prelim}

\begin{lem}\label{lem1}  Let $K\unlhd H \leq T\leq G.$ Assume the triple $(G,H,K)$ is special in $G.$ \\
(a) The  triple $(T,H,K)$ is special in $T.$\\
(b)  If $L/K^G\unlhd HK^G/K^G$ and triple $(G,H,L\cap H)$ is special in $G$ then the   triple $(G/K^G,HK^G/K^G,L/K^G)$ is special in $G/K^G.$ In particular, if $K\unlhd G$ then the triple $(G/K,H/K,L/K)$ is special in $G/K.$\\
(c) If $H$ is an $NR$-subgroup of $G$ then $HK^G/K^G$ is an $NR$-subgroup of
$G/K^G.$\\
(d)  If $K\unlhd G$ and every non-nilpotent maximal subgroup of $G$ is either normal or $NR$ in $G,$ then every non-nilpotent maximal subgroup of $G/K$ is either normal or $NR$ in $G/K.$\\
$(e)$ Let $K\unlhd H \leq G.$ If there is a subgroup
$L$ such that $G=HL$ and $H\cap L=K$ then the triple $(G,H,K)$ is
special.
\end{lem}
\begin{proof} $(a)-(c)$ are in Lemma $4$ in Berkovich\cite{ber}. As nilpotence and self-normalizing are preserved under taking quotient group, this together with $(c)$ yield $(d).$ Finally $(e)$ is Lemma $9$ in Berkovich\cite{ber}.
\end{proof}

\begin{lem}\label{lem2} \emph{(Tate\cite{tate}).} Let $H\unlhd G$ and $P\in Syl_p(G).$ If $H\cap P\leq \Phi(P)$ then $H$ is $p$-nilpotent.
\end{lem}

\begin{thm}\label{zsi}\emph{(Zsigmondy\cite{Zsi}).} Let $q$ and $n$ be integers with $q\geq 2$ and $n\geq 3.$ If $(q,n)\neq (2,6)$ then there is a prime $r$ such that $r\mid q^n-1$ but $r$ does not divide $q^i-1$ for $i<n.$
\end{thm}
We call such a prime $r$ a \emph{primitive prime divisor} of $q^n-1$ and denoted by $q_n.$
\begin{lem}\label{lem3} Let $t\geq 0$ be an integer.\\
$(a)$ $2^t-1,$ and $ 2^{t-1}-1$ are both prime powers if and only if  $t=3.$\\
$(b)$ $2^t+1,$ and $ 2^{t+1}+1$ are both prime powers if and only if  $0\leq t\leq 3.$
\end{lem}
\begin{proof} $(a)$ Observe first that $t$ and $t-1$ are both primes if and only if $t=3,$ in which case $2^t-1=7$ and $2^{t-1}-1=3$ are both primes. If $t=6$ then $2^6-1=9.7$ is not a prime power. Assume that $t>3$ and $t\neq 6.$ Then either $t$ or $t-1$ is not prime. Assume that $n>3, n\neq 6$ is not a prime. We will show that $2^n-1$ is not a prime power. By way of contradiction, assume that $2^n-1=p^m$ for some $m.$  By Theorem \ref{zsi}, $p$ is a primitive prime divisor of $2^n-1.$
Let $s<n$ be a non-trivial prime divisor of $n$ and write $n=sa.$ Then $p^m=2^n-1=2^{sa}-1$ is divisible by $2^s-1.$ It follows that $p$ is a divisor of $2^s-1$ with $s<n,$ contradicting to the definition of primitive prime divisor. Thus if either $t$ or $t-1$ is not prime then either $2^t-1$ or $2^{t-1}-1$ is not a prime power.

$(b)$ If $0\leq t\leq 3,$ then we can check that both $2^t+1$ and $2^{t+1}+1$ are prime powers. Assume that $t\geq 4,$ and $2^t+1, 2^{t+1}+1$ are both prime powers. We see that $t$ or $t+1$ must be odd. Now let $n$ be an odd integer with $n\geq 4.$ We will show that $2^n+1$ is not a prime power. By way of contradiction, assume that $2^n+1=p^m$ for some prime $p.$ Since $n$ is odd, we have $3|2^n+1$ so that $p=3.$ Thus $2^n+1=3^m$ or equivalently $3^m-1=2^n.$  As $n\geq 4,$ we see that $m\geq 3.$ By Zsigmondy's Theorem, $3^m-1$ has a primitive prime divisor which is not a divisor of  $3-1=2,$ a contradiction. Thus $2^t+1$ and $2^{t+1}+1$ cannot be both prime powers when $t\geq 4.$
\end{proof}

\begin{lem}\label{lem4}\emph{(Thompson \cite[Bemerkung II$.7.5$]{hu1}).} If $G$ is a minimal simple group then $G$ is isomorphic to one of the
following groups:\\
$(1)$ $L_2(p),$ $p>3$ is a prime, and $p^2-1\not\equiv 0$ \emph{(mod $5$)};\\
$(2)$ $L_2(3^r),$ $r$ is an odd prime;\\
$(3)$ $L_2(2^r),$ $r$ is a prime;\\
$(4)$ $Sz(2^r),$ $r$ is an odd prime;\\
$(5)$ $L_3(3).$
\end{lem}

\begin{lem}\label{lem5} $(i)$ \emph{(Corollary $2.2,$ King\cite{king}).} Assume that $G\cong L_2(q)$ is a minimal simple group.
If $M$ is a maximal subgroup of $G$ then $M$ is one of the following groups:\\
$(a)$ $D_{q-1}$ for $q\geq 13$ odd and $D_{2(q-1)}$ for $q$ even;\\
$(b)$ $D_{q+1}$ for $q\neq 7,9$ odd and $D_{2(q+1)}$ for $q$ even;\\
$(c)$ a Frobenius group of order $q(q-1)/2$ for $q$ odd and $q(q-1)$ for $q$ even;\\
$(d)$ $S_4$ when $q\equiv \pm 1$ \emph{(mod $8$)} with $q$ prime or $q=p^2$ and $p\equiv \pm 3$ (mod $8$);\\
$(e)$ $A_4$ when $q\equiv \pm 3$ \emph{(mod $8$)} with $q$ prime;\\
$(ii)$ \emph{(Theorem $9,$ Suzuki\cite{su2}).} Assume that $G\cong Sz(q),q=2^r,r$ an odd prime. The maximal subgroup of $G$ are as follows:\\
$(a)$ a Frobenius group of order $q^2(q-1);$\\
$(b)$ a dihedral subgroup of order $2(q-1);$\\
$(c)$ a Frobenius group of order $4(q\pm t+1),$ with $t^2=2q.$\\
$(iii)$ \emph{(Conway et al.\cite{atlas}, page 13).}
Assume that $G\cong L_3(3).$  The maximal subgroup of $G$ are as follows:\\
$(a)$ a  group of order $3^2:2S_4$\\
$(b)$ a group of order $13:3$\\
$(c)$ $S_4.$
\end{lem}

Let $p$ be a prime. Recall that $O^p(G)$ is the smallest normal subgroup of $G$ such that $G/O^p(G)$ is a $p$-group, or equivalently, $O^p(G)$ is a subgroup of $G$ generated by all $p'$-elements in $G.$ Also $F(G)$ is the \emph{Fitting subgroup} of $G,$ that is the largest nilpotent normal subgroup of $G.$

\begin{thm}\label{BB} \emph{(Satz Baumann\cite{bau}).} Let $G$ be a non-solvable group possessing a nilpotent maximal subgroup. Then $O^2(G/F(G))$ is a direct product of simple groups with dihedral $2$-Sylow subgroups. The simple groups that can appear are $L_2(q)$ with $q=9$ or $q$ a prime of the form $2^m\pm 1>5.$
\end{thm}

\begin{remark} We have $L_2(9)\cong A_6,$ $Out(A_6)\cong \Z_2\times \Z_2$ and $S_6\cong A_6.2_1.$ The $2$-Sylow subgroups of $A_6$ and $S_6$ are not maximal but the $2$-Sylow subgroups of other extensions of $A_6$ are maximal. \emph{(see Conway et al.\cite{atlas} page $4$).}
\end{remark}

For any $p$-group $P,$ we denote by $\mathcal{A}(P)$ the set of abelian subgroups of $P$ of maximal order. The \emph{Thompson subgroup} $J(P)$ is a subgroup of $P$ generated by all members of $\mathcal{A}(P).$ It is well known that $J(P)$ is characteristic in $P,$ $Z(J(P))$ is characteristic in $J(P)$ and hence $Z(J(P))$ is characteristic in $P.$ The following theorem gives a sufficient condition for the existence of a normal $p$-complement in finite groups.
\begin{thm}\label{GT} \emph{(Glauberman-Thompson, Theorem $8.3.1,$ Gorenstein\cite{gor}).} Let $P\in Syl_p(G)$ with $p$ odd. If $N_G(Z(J(P)))$ has a normal $p$-complement then $G$ has a normal $p$-complement.
\end{thm}

Recall that a group $G$ is called \emph{minimal non-nilpotent} if $G$ is non-nilpotent but every maximal subgroup of $G$ is nilpotent. The structure of minimal non-nilpotent groups are given by the following theorem:
\begin{thm}\label{Sch} \emph{(O.J. Schmidt, $9.1.9,$ Robinson\cite{rob}).} Assume that $G$ is a minimal non-nilpotent group. Then\\
$(a)$ $G$ is solvable;\\
$(b)$ $|G|=p^mq^n$ where $p,q$ are distinct primes. Morever, there is a unique $p$-Sylow subgroup $P$ and a cyclic $q$-Sylow subgroup $Q.$ Hence $G=PQ$ and $P\unlhd G.$
\end{thm}

\begin{lem}\label{lem6} $(a)$ If $S$ is a non-abelian simple group and $S\unlhd A\leq Aut(S),$ then there exists a subgroup $K\leq S$ such that $A=SN_A(K)$ and every proper-over group of $K$ in $S$ is local in $S.$ Moreover when $S$ is a finite group of Lie type and $S\neq {}^2F_4(2)',$ we can choose $K$ to be a Borel subgroup of $S.$\\
$(b)$ If every maximal subgroup of $G$ is an $NR$-subgroup then $G$ is solvable.
\end{lem}
\begin{proof} $(a)$ follows from Theorem $1.2$\cite{tong} and its proof. $(b)$ is Theorem $1.1$\cite{tong}.
\end{proof}

The following results are obvious.
\begin{lem}\label{lem7} Let $A$ be a normal subgroup of prime order of $G.$ Then $G$ is supersolvable if and only if $G/A$ is supersolvable.
\end{lem}

\begin{lem}\label{lem8} Let $H$ be a normal subgroup of $G.$ Assume that $H$ is $p$-nilpotent with a normal $p$-complement $K.$ Then $K\unlhd G.$
\end{lem}
\begin{lem}\label{lem9} Let $t\geq 0$ be an integer.\\
$(a)$ $D_{2t}$ is nilpotent if and only if $t$ is a power of $2.$\\
$(b)$ $D_{2t}$ is minimal non-nilpotent if and only if $t$ is an odd prime.
\end{lem}
 \section{$NR$-subgroups and solvability}
 \label{sol}

The following result is a generalization of exercise $9.1.10$ in Robinson\cite{rob}.
\begin{lem}\label{nr1} If every self-normalizing maximal subgroup of $G$ is nilpotent then $G$ is solvable.
\end{lem}
\begin{proof} As nilpotent groups are solvable, we can assume that $G$ is not nilpotent. If every maximal subgroup of $G$ is self-normalizing then every maximal subgroup of $G$ is nilpotent by hypothesis, hence  $G$ is a minimal non-nilpotent group, and so by Theorem \ref{Sch}, $G$ is solvable. Thus $G$ contains a maximal subgroup which is not self-normalizing, hence  $G$ is not simple. If $A$ is any non-trivial normal subgroup of $G$ then $G/A$ satisfies the hypothesis of the lemma, so that by induction, $G/A$ is solvable. Therefore $G$ has a unique minimal normal subgroup $N$ and $G/N$ is solvable. Assume $N$ is not solvable, hence $N$ is a direct product of some non-abelian simple groups. Let $P\in Syl_p(N),$ with $p$ odd. By Frattini's argument, we have $G=N_G(P)N.$ Clearly $N_G(P)<G,$ hence there exists a maximal subgroup $M$ of $G$ such that $N_G(P)\leq M< G.$ We have $G=N_G(P)N=MN$ and hence $M$ is not normal in $G.$ Therefore $M$ is nilpotent. Let $Q\in Syl_p(M).$ As $M$ is nilpotent, we have $Q\unlhd M$ so that $M=N_G(Q).$ We will show that $Q\in Syl_p(G).$ By way of contradiction, assume that $Q$ is not a $p$-Sylow subgroup of $G.$ Let $S$ be a $p$-Sylow subgroup of $G$ containing $Q.$ We have $Q<N_S(Q),$ and hence $Q<N_S(Q)\leq N_G(Q)=M,$ a contradiction as $Q\in Syl_p(M).$ As $Q\unlhd M,$ the Thompson subgroup $J(Q)$ is characteristic in $Q$ and $Z(J(Q))$ is characteristic in $J(Q),$ and so $Z(J(Q))$ is characteristic in $Q.$ It follows that $Z(J(Q))\unlhd M$ so that $M=N_G(Z(J(Q))).$ Since $M$ is nilpotent, it is $p$-nilpotent and hence by Theorem \ref{GT}, $G$ has a normal $p$-complement $H.$ By the uniqueness of $N,$ we have $N\leq H.$ This is a contradiction as $p$ divides the order of $N$ but not that of $H.$ Thus $N$ is solvable and so $G$ is solvable.
\end{proof}

{\bf Proof of Theorem \ref{th1}}.
Let $G$ be a minimal counter-example to Theorem \ref{th1}.
Assume first that $G$ is non-abelian simple. Then every maximal subgroup of $G$ is either nilpotent or $NR.$ By Lemma \ref{lem6}$(b),$ $G$ contains a nilpotent maximal subgroup $H.$ By Theorem \ref{BB}, we have $G\cong L_2(q)$ for $q$ a prime of the form $2^m\pm 1>5.$ By Lemma \ref{lem5}$(i),$ $G$ has a maximal subgroup $B$ of order $q(q-1)/2$ which is a Frobenius group with Frobenius kernel $U$ of order $q$ and a Frobenius complement $T$ with $|T|=(q-1)/2>1.$ Since $B$ is not nilpotent, from hypothesis, $B$ is an $NR$-subgroup of $G.$ As $G$ is simple and $U$ is non-trivial, we have $G=U^G,$ hence $U^G\cap B=G\cap B=B>U,$ contradicting to the fact that $B$ is an $NR$-subgroup of $G.$ Therefore $G$ is not non-abelian simple.

Since the hypothesis of the Theorem inherits to proper quotient of $G,$ $G$ has a unique minimal normal subgroup $N$ which is not solvable. Then $$N=S^{x_1}\times S^{x_2}\times\cdots \times S^{x_t},$$ where $S$ is a non-abelian simple group, and
$x_1, x_2,\cdots,x_t\in G.$ Let $K$ be the subgroup of $S$ obtained from Lemma  \ref{lem6}$(a)$, and let $$R=K_1\times K_2\times\cdots
\times K_t,$$ where $K_i=K^{x_i},i=1,\cdots,t.$  Then $R$ is a non-trivial proper subgroup of $N.$ Since $N$ is a unique minimal normal subgroup of $G,$ $N_G(R)<G.$ We will show that $G=N_G(R)N.$ For any $g\in G,$ since $N^g=N,$ there exists a permutation $\pi\in S_t$ such that $S^{x_ig}=S^{x_{i\pi}}.$ Let $g_i=x_i g x_{i\pi}^{-1}.$ Then $g_i\in N_G(S).$ We have $$R^g=K^{x_1g}\times K^{x_2g}\times\cdots \times K^{x_tg}=K^{g_1x_{1\pi}}\times K^{g_2x_{2\pi}}\times\cdots \times K^{g_tx_{t\pi}}=$$ $$=K^{g_{1\pi^{-1}}x_1}\times
K^{g_{2\pi^{-1}}x_2}\times\cdots \times K^{g_{t\pi^{-1}}x_t}=K^{h_1x_1}\times K^{h_2x_2}\times\cdots\times K^{h_{t}x_t}=$$ $$=K^{x_1s_1}\times K^{x_2s_2}\times\cdots
\times K^{x_ts_t}=K_1^{s_1}\times K_2^{s_2}\times\cdots\times K_t^{s_t},$$ where $K^{g_{i\pi^{-1}}}=K^{h_i}$ with $h_i\in S$ by Lemma \ref{lem6}$(a)$,  and $s_i=h_i^{x_i}\in S^{x_i}.$ Let $s=s_1.s_2\dots s_t\in N.$ Since $[S^{x_i},S^{x_j}]=1$ if $i\not=j\in\{1,2,\dots, t\},$ $K_i^s=K_i^{s_i}.$ Thus $R^g=R^s, $ where $s\in N.$ Therefore $G=N_G(R)N.$ Let $M$ be a maximal subgroup of $G$ containing $N_G(R).$ We have $G=MN.$ 

Clearly $M$ is not normal in $G,$ otherwise $N\leq M$ and hence $G=M,$ a contradiction.

Assume that $M$ is nilpotent. As $F(G)=1,$ by Theorem \ref{BB}, $O^2(G)\unlhd G$ is a product of simple groups, it follows that $N\leq O^2(G)$ and hence $S\cong L_2(q)$ for $q=9$ or $q$ a prime of the form $q=2^m\pm 1>5.$  We will show that  $\pi(M)\cap \pi(N)=\{2\}.$ Assume $p\in \pi(M)\cap \pi(N)$ and $p$ is odd. As in proof of Lemma \ref{nr1}, the $p$-Sylow subgroup $P$ of $M$ is also a $p$-Sylow subgroup of $G.$ Apply Glauberman-Thompson Theorem, $G$ has a normal $p$-complement. This leads to a contradiction as in the proof of the previous lemma. Thus $M\cap N$ is a $2$-group. Since $R\leq M\cap N,$ $R$ must be a $2$-group and hence $K$ is a $2$-group. However by the choice of $K$ in Lemma \ref{lem6}$(a),$ $|K|=q(q-1)/2$ which is not a $2$-group as $q$ is odd. Thus $M$ is not nilpotent.

Therefore we can assume that $M$ is not nilpotent nor normal in $G.$ Then $M$ is an $NR$-subgroup of $G.$ The argument below is exactly the same as in the last part of the proof of Theorem $1.1$ in \cite{tong}.
Let $Q=M\cap N.$ We have $G=MN,$ and $Q=M\cap N\unlhd M.$ As $G/N=MN/N\simeq M/Q,$
$M/Q$ is solvable. If $Q$ is solvable then $M$ is solvable. By Proposition $7$ in \cite{ber}, $G=M,$ a contradiction. Thus $Q$ is non-solvable. Let $L$ be any non-trivial normal subgroup of $M.$ Since $M$ is an $NR$-subgroup  of $G,$ we have $L=L^G\cap M.$ It follows from the fact that $N$ is the unique minimal normal subgroup of $G,$ $N\leq L^G.$ We have $Q=N\cap M\leq L^G\cap M=L.$ We conclude that $Q$ is a minimal normal subgroup of $M.$ Since $Q$ is a minimal normal subgroup of $M$ and $Q$ is non-solvable, $Q=W_1\times W_2\times \cdots\times W_k,$ where $W_i\simeq W$ for all $1\leq i\leq k$ and $W$ is a non-abelian simple group. Suppose that there exists $j\in \{1,2,\dots, t\}$ such that $S^{x_j}\leq Q.$ As $S^{x_j}$ is normal in $N,$ $(S^{x_j})^G=(S^{x_j})^{NM}=(S^{x_j})^M\leq M.$ However  as $(S^{x_j})^G=N,$ $G=MN=M,$ a contradiction. Therefore $S^{x_j}\cap Q<S^{x_j}$ for any $j\in\{1,2,\dots, t\}.$ Since $K_j\leq S^{x_j}\unlhd N,$ we have $K_j\leq S^{x_j}\cap Q\unlhd Q.$ As $Q$ is a direct product of non-abelian simple groups and $S^{x_j}\cap Q$ is a non-trivial normal subgroup of $Q,$ there exists a non-empty set $J\subseteq\{1,2,\dots, t\}$ such that $S^{x_j}\cap Q=\prod_{i\in J}W_i.$ Hence $K_j\leq \prod_{i\in J}W_i<S^{x_j},$ and so $K\leq \prod_{i\in J}W_i^{x_j^{-1}}<S,$ where $W_i^{x_j^{-1}}$ are non-abelian simple for any $i\in J.$ However by Lemma \ref{lem6}$(a),$ $\prod_{i\in J}W_i^{x_j^{-1}}$ is local in $S.$ This final contradiction completes the proof.   $\hfill \square$

{\bf Proof of Theorem \ref{th2}}.
Let $M$ be any maximal subgroup of $G.$ From hypothesis, every non-nilpotent maximal subgroup $H$ of $M$ is either subnormal or $NR$ in $G.$ By Lemma \ref{lem1}$(a)$ and the maximality of $H$ in $M,$ $H$ is either normal or $NR$ in $M,$ so that $M$ satisfies the hypothesis of Theorem \ref{th1}. Thus $M$ is solvable. It follows that every maximal subgroup of $G$ is solvable. If $N$ is any non-trivial proper normal subgroup of $G,$ then $N$ is solvable and hence $N\leq S(G).$ By Lemma \ref{lem1}$(d),$ $G/N$ satisfies the hypothesis, so that $(G/N)/S(G/N)$ is trivial or isomorphic to $A_5.$
Since $S(G/N)=S(G)/N,$ we have $(G/N)/S(G/N)\cong (G/N)/(S(G)/N)\cong G/S(G).$

Thus $G/S(G)$ is either trivial or isomorphic to $A_5,$ and we are done. Therefore we can assume that $G$ is simple and so $G$ is a minimal simple group. By Theorem \ref{Sch}$(a),$ $G$ contains a  non-nilpotent maximal subgroup. Now let $M$ be any non-nilpotent maximal subgroup of $G.$  Let $H$ be any maximal subgroup of $M,$ as $H$ is not subnormal in $G,$ $H$ is either nilpotent or $NR$ in $G.$ We will show that if $H$ is $NR$ in $G$ then $H$ is of prime order. Assume that $H$ is $NR$ in $G.$ Let $A$ be any non-trivial normal subgroup of $H.$ As $H$ is an $NR$-subgroup of $G,$ the triple $(G,H,A)$ is special in $G$ so that $A^G\cap H=A.$ As $G$ is simple, we have $A^G=G,$ hence $H=A,$ so that $H$ is a simple group. Since $H$ is solvable, $H$ must be cyclic of prime order. Therefore $M$ is a minimal non-nilpotent group. By Theorem \ref{Sch}$(b),$ $|M|=p^mq^n$ where  $p,q$ are distinct primes and $m,n$ are a positive integer, the $p$-Sylow subgroup $P$ of $M$ is normal in $M$ while the $q$-Sylow subgroup $Q$ of $M$ is cyclic. Now by Lemma \ref{lem4}, we consider the following cases:

$(a)$ Case $G\cong L_2(q),q>3$ prime, and  $q^2-1\not\equiv 0$ (mod $5$) or $G\cong L_2(q),q=3^r,r$ odd prime. If $q=5$ then $G=L_2(5)\cong A_5$ and we are done. If $q=7$ then $G=L_2(7)$ and since $7\equiv -1(\mbox{mod } $8$),$ $S_4$ is a maximal subgroup of $G.$ However $S_4$ is neither nilpotent nor minimal non-nilpotent as it contains a non-nilpotent subgroup $S_3.$ Therefore we can assume that $q\geq 13.$ By Lemma \ref{lem5}$(i),$ the Frobenius group $F$ of order $q(q-1)/2$ and the dihedral groups $D_{q-1}, D_{q+1}$ are maximal subgroup of $G.$ Since $q\geq 13$  and $1<(q-1)/2<q,$ $F$ is not nilpotent, and so $(q-1)/2=s^a, s$ prime. We have $q+1=2(s^a+1).$ If $D_{q+1}$ is nilpotent then $q+1=2(s^a+1)=2^t$ for some integer $t\geq 0,$ by Lemma \ref{lem9}$(a).$ Hence $q=2^t-1$ and $(q-1)/2=2^{t-1}-1$ are both prime powers. Since $q\geq 13,$ we have $t\geq 4$ and so $2^t-1$ and $2^{t-1}-1$ cannot be prime powers at the same time by Lemma \ref{lem3}$(a)$. Thus $D_{q+1}$ is not nilpotent, and hence  $D_{q+1}$ is minimal non-nilpotent. Thus $(q+1)/2=s^a+1$ must be an odd prime by Lemma \ref{lem9}$(b).$ It follows that $s$ is an even prime and so $s=2.$ Hence $q=2^{a+1}+1$ and $2^a+1$ are both prime powers. By Lemma \ref{lem3}$(b),$ we have $a=0,1,2,3$ and so $q=3,5,9,17$ respectively. As $q\geq 13,$ we have $q=17.$ However $(17+1)/2=3^2$ is not an odd prime. Thus these cases cannot happen unless $q=5.$

$(b)$ Case $G\cong L_2(2^r),r$ prime. As $L_2(4)\cong L_2(5)\cong A_5,$ we can assume that $r$ is an odd prime. By Lemma \ref{lem5}$(i),$ $G$ has a maximal subgroup isomorphic to $D_{2(q+1)}.$ Clearly $D_{2(q+1)}$ is not nilpotent and hence it must be minimal non-nilpotent. By Lemma \ref{lem9}$(b),$ $q+1=2^r+1$ is prime. As $r$ is odd, $2^r+1$ is divisible by $3.$ Thus $2^r+1=3$ so that $r=1,$ a contradiction.

$(c)$ Case $G\cong Sz(q),q=2^r,r$ odd prime. By Lemma \ref{lem5}$(ii),$ and Theorem $3.10$\cite{hub}, $G$ has a maximal subgroup $M=N_G(A),$ where $A$ is cyclic of order $q+t+1$ with $2t^2=q,$ $|M:A|=4$ and $M$ is a Frobenius group with Frobenius kernel $A.$ Moreover $C_G(u)=A$ for any $1\neq u\in A.$ As $q+t+1$ is odd and $C_G(u)=A$ for any $1\neq u\in A,$ we see that $M$ is not nilpotent. As $|M:A|=4,$ $M$ contains a subgroup of index $2$ which is not nilpotent. Thus $M$ is not a minimal non-nilpotent.

$(d)$ Case $G\cong L_3(3).$   By Lemma \ref{lem5}$(iii),$ $G$ has a maximal subgroup $M$ which is isomorphic to $S_4.$ As in case $(a),$  $S_4$ is neither nilpotent nor minimal non-nilpotent.

The proof is now complete.$\hfill\square$

{\bf Proof of Corollary \ref{cor}} As $G$ satisfies the hypothesis of Theorem \ref{th2}, $G$ is solvable or $G/S(G)\cong A_5.$ If the first possibility holds then we are done. Thus we assume that $G/S(G)\cong A_5.$ As the hypothesis of the corollary inherits to quotient group and since $A_5$ is simple, it follows that every $2$-maximal subgroup of $A_5$ is an $NR$-subgroup. As $V_4\unlhd A_4\leq A_5$ and $V_4$ is maximal in $A_4,$ $A_4$ is maximal in $A_5,$ we deduce that $V_4$ is a $2$-maximal subgroup of $A_5.$ Let $K$ be any subgroup of order $2$ in $V_4.$ Then $K\unlhd V_4\leq A_5.$ It follows that $K^{A_5}\cap V_4=K.$ However as $A_5$ is simple, and $1\neq K^{A_5}\unlhd A_5,$ we have $K^{A_5}=A_5,$ and hence $K^{A_5}\cap V_4=V_4>K.$ This contradicts to the fact that $H$ is an $NR$-subgroup of $A_5.$ Therefore $G$ is solvable.
The proof is now complete.$\hfill\square$

\section{Normal Complement}
\label{norm}

{\bf Proof of Proposition \ref{nc1} } Let $K=P^G\unlhd G.$ As  $(G,N,P)$ is special in $G,$ $K\cap N=P.$   Also $N_K(P)=K\cap N=P,$ $P$ is a self-normalizing $p$-subgroup of  $K$ so that $P$ is a self-normalizing $p$-Sylow subgroup of $K.$ By Frattini's argument, we have $G=NK.$ Let $L=\Phi(P)^G.$ We have $L\unlhd G$ and since the triple $(G, N,\Phi(P))$ is special in $G,$ $\Phi(P)\leq P\cap L\leq N\cap L=\Phi(P),$ hence $P\cap L=\Phi(P).$  As $L\leq K$ and $PL/L\cong P/(P\cap L)\cong P/\Phi(P),$ $PL/L$ is an elementary abelian $p$-Sylow subgroup of $K/L.$ We will show that $PL/L$ is self-normalizing in $K/L.$ In fact, assume $Lg\in K/L$ normalizes $PL/L$ where $g\in K.$ Then we have $P^g\leq PL\leq K.$ Since $P$ is a $p$-Sylow subgroup of $K,$ by Sylow's Theorem, $P^g=P^u$ for some $u\in PL.$ Indeed, we can take $u\in L.$ Thus $g\in N_K(P)L=PL.$ This proves our claim. Now by Burnside Normal $p$-Complement Theorem, (Theorem $7.4.3$ Gorenstein\cite{gor}), $PL/L$ has a normal $p$-complement $H/L$ in $K/L.$ It follows that $K=PLH=PH$ and $PL\cap H=L.$ As $K/L\unlhd G/L$ and $H/L$ is a normal $p$-complement in $K/L,$ by Lemma \ref{lem8}, $H/L\unlhd G/L$ and hence $H\unlhd G.$ Observe that $$P\cap H\leq P\cap (PL\cap H)=P\cap L\leq N\cap L=\Phi(P)$$ as the triple $(G,N,\Phi(P))$ is special in $G.$ By Lemma \ref{lem2}, $H$ has a normal $p$-complement $T$ in $H.$ Thus $K=PH=PT, P\cap T=1$ and so $$G=NK=NPT=NT$$ and $$N\cap T=(N\cap PT)\cap T=(N\cap K)\cap T=P\cap T=1.$$ Applying Lemma \ref{lem8} again for $H\unlhd G,$ we have $T\unlhd G.$ Thus $T$  is a normal complement for $N$ in $G.$ The proof is now complete. $\hfill\square$

{\bf Proof of Theorem \ref{th4}}. By Proposition  \ref{nc1}, $G=NL$ for some normal $p'$-subgroup $L$ of $G.$ Assume that $N$ is $p$-nilpotent. Then $N=S.K$ where $P\leq S\in Syl_p(N)$ and $K$ is a normal $p'$-subgroup of $N.$  We have $G=NL=SKL.$ Clearly $KL\unlhd G,$ $KL$ is a $p'$-subgroup and so $S\in Syl_p(G).$ Hence $G$ is $p$-nilpotent.$\hfill\square$

{\bf Proof of Theorem \ref{th5}}. Assume that $H$ is supersolvable. Let $A$ be a minimal normal subgroup of $H.$ Then $A$ is a cyclic subgroup of prime order $p.$ Assume first that $A\unlhd G.$ Then $H/A$ is an $NR$-subgroup of $G/A$ by Lemma \ref{lem1}$(c).$ Moreover $H/A$ is supersolvable and has prime order in $G/A.$ By induction we have $G/A$ is supersolvable. The result follows from Lemma \ref{lem7}. Now assume that $A$ is not normal in $G.$ We have $H=N_G(A).$ By Proposition \ref{nc1},  $H$ has a normal complement $K$ in $G$ so that $G=HK$ and $H\cap K=1.$ It follows that $|K|=|G:H|$ is a prime. Applying Lemma \ref{lem7} again, $G$ is supersolvable. $\hfill\square$

\subsection*{Acknowledgment}
I am grateful to Professor Christopher W. Parker for drawing my attention to the paper of B. Baumann \cite{bau}.

\end{document}